\documentclass{lmsedition}
\usepackage{blindtext}

\makeatletter
\def\fakeH@#1#2{%
  \vtop{\m@th\ialign{##\cr
    \hfil$#1\kern-.8pt\operator@font \textsf{o}$\hfil\cr
    \noalign{\nointerlineskip\kern-7\ex@}#2\cr
    \noalign{\nointerlineskip\kern-\ex@}\cr}}%
}
\def\fakHo{%
  \mathop{\mathpalette\fakeH@{\'{}\kern-.6pt\'{}}}\nmlimits@
\!}

\def\set@curr@file#1{%
	\begingroup
	\escapechar\m@ne
	\xdef\@curr@file{\expandafter\string\csname #1\endcsname}%
	\endgroup
}
\def\quote@name#1{"\quote@@name#1\@gobble""}
\def\quote@@name#1"{#1\quote@@name}
\def\unquote@name#1{\quote@@name#1\@gobble"}
\makeatother

\makeatother

\usepackage[utf8]{inputenc}
\usepackage[english]{babel}
\usepackage{mathtools, bm}
\usepackage{here}
\usepackage{tikz}
\usetikzlibrary{matrix,arrows,automata}

\newcommand\Kbarz{\overline{\mathbb{K}}(z)}

\newcommand\Kbar{\overline{\mathbb{K}}}
\newcommand\Qbar{\overline{\mathbb{Q}}}
\newcommand\fiz{f_{1}(z),\ldots, f_{n}(z)}
\newcommand\fialpha{f_{1}(\alpha),\ldots, f_{n}(\alpha)}
\newcommand\Xis{X_{1},\ldots, X_{n}}

\newcommand\FqTc{\mathbb{F}_{q}(T)}
\newcommand\FqTa{\mathbb{F}_{q}[T]}

\begin{document}

\articletitle{Historical developments of the theory of transcendence and algebraic independence}

\articleauthor{Gwladys Fernandes}

\begin{articleabs}
\textit{Qu'est la majesté de ce qui finit, auprès des départs titubants, des désordres de l'aurore ?}\footnote{\textit{Flore et Pomone}, Colette. Translation:\textit{What is the majesty of what ends, compared to the staggering departures, to the disorders of the dawn?}}. 
Let this question from the french writer Colette, guide us to walk through the promising beginnings, fascinating evolutions and results, challenging open problems and perspectives of an exciting theory that is far from revealing all its secrets.
\end{articleabs}

\begin{twoblock}


\section*{Historical methods}

The notion of a transcendental number has been slow to construct. This concept reflects the discoveries of confusing numbers, which have escaped to the control and the understanding of mathematicians for centuries. The term \textit{transcendante} is finally first introduced by G. Leibniz in 1682 to refer to certain curves, as graphs of sinus and cosinus, in contrast to {algebraic} curves. But it was not until 1844 that J. Liouville gave a formal definition of these puzzling numbers. These are the complex numbers which are not roots of any polynomial with integers as coefficients. The complement of this set is the one of algebraic numbers over $\mathbb{Q}$, denoted as $\overline{\mathbb{Q}}$. 

\textbf{From J. Liouville to K. Weierstrass}

Among the first examples of transcendental numbers provided by J. Liouville, there is $\sum_{k=1}^{+\infty}2^{-k!}$. The proof relies on a fundamental theorem of diophantine approximation from the author, called Liouville inequality. Let $\alpha\in\Qbar$, and let $n$ be its degree, that is the minimal degree among the ones of the polynomials $P(X)\in\mathbb{Q}[X]$ such that $P(\alpha)=0$. Then, the inequality of Liouville states that there exists a constant $c>0$ such that for all $p\in\mathbb{Z},q\in\mathbb{N}^{*}$ with $p/q\neq\alpha$, we have
\begin{equation}
\label{Liouville}
\left|\alpha-\frac{p}{q}\right|>\frac{c}{q^{n}}.
\end{equation}   
Indeed, the decrease of $2^{-k!}$ makes the series \textit{too much well approached} by rational numbers to be algebraic over $\Qbar$. 

Note that Liouville inequality is improved in 1955 by K. Roth who replaces $n$ by $2+\epsilon$ in \eqref{Liouville}, for any $\epsilon>0$. Furthermore, it is generalised in 1972 as a simultaneous approximation result to several algebraic numbers by rational numbers with the same denominator by W. Schmidt and his sub-space theorem.  

Then, in 1874, G. Cantor establishes that the set of real numbers is not countable, whereas the set $\Qbar$ of algebraic numbers is. Hence, he shows that \textit{almost all} real or complex numbers are transcendental. 

But, one of the first big events of the theory of transcendence already took place one year before with the proof of the transcendence over $\Qbar$ of the number $e$ by C. Hermite. The essential ingredient, nowadays called \textit{Padé approximant}, is the key of the evolution of the theory of transcendence. Following the same method, F. Lindemann proves in 1882 that $\pi$ is transcendental over $\Qbar$. Then, the author establishes the transcendence of the values of the exponential function at non-zero algebraic points. F. Lindemann also states a more general result, but without a complete proof, which is written later by K. Weierstrass. 

This result, called the Lindemann-Weierstrass's theorem states that, given linearly independent algebraic numbers $\rho_{1},\ldots,\rho_{n}$ over $\mathbb{Q}$, the numbers  
$e^{\rho_{1}},\ldots,e^{\rho_{n}}$ are algebraically independent over $\Qbar$. This means that there is no non-zero polynomial $P(\Xis)\in\Qbar[\Xis]$ such that $P(e^{\rho_{1}},\ldots,e^{\rho_{n}})=0$.

The Lindemann-Weierstrass's theorem provides an example of algebraically independent functions that take algebraic independent values at non-zero algebraic numbers. This transfer of algebraic independence from functions to their values at some algebraic points is not true for arbitrary functions. However, there exist at least three cases in which this holds : differential systems, Mahler systems and $\tau$-difference systems (in positive characteristic). In what follows, we will present these three settings, starting with the ones which appear in characteristic zero. Our aim is to stress the major developments they bring to the theory of transcendence and algebraic independence.

\begin{bandbox}[dcolor]
		\section*{Padé approximants}
	In Hermite's proof, Padé approximants refer to the simultaneous approximation to integral values of the exponential function by rational fractions with the same denominator. In the proof of the Lindemann-Weierstrass's theorem, they arise from functional linear approximating forms for $e^{\rho_{1}z},\ldots,e^{\rho_{n}z}$. This latter construction is possible thanks to the famous Pigeon hole principle, which gives rise to some remarkable results. The first ones are due to C. Siegel, in 1929. For example, the author proves that the set of algebraic integral points on an algebraic curve of genus $g>0$ over $\mathbb{Z}$ is finite. Note that the second crucial argument of his proof is a refinement of the Liouville inequality. 
\end{bandbox}

\textbf{C. Siegel work and A. Shidlovskii contribution}

In the same paper of 1929 mentioned above, C. Siegel introduces a new class of functions, called $E$-functions, that are in particular solutions of linear differential systems with coefficients in $\Qbar(z)$. Examples of $E$-functions are the exponential and the Bessel functions, along with some hypergeometric functions. The work of C. Siegel is completed later by A. Shidlovskii and gives the following fundamental result, known as the Siegel-Shidlovskii's theorem. Let $(f_{1}(z), \ldots, f_{n}(z))$ be a vector of $E$-functions, which satisfies a linear differential system. Let $\alpha$ be a non-zero algebraic number that is not a pole of the matrix of the system. Then, the transcendence degree of $\{f_{1}(\alpha), \ldots, f_{n}(\alpha)\}$ over $\Qbar$ is equal to the one of $\{f_{1}(z), \ldots, f_{n}(z)\}$ over $\Qbar(z)$. Recall that the transcendence degree of a family is the maximal number of algebraically independent numbers it contains. Recently, Y. André recovers this statement by developing a new general Galois theory based on affine quasi-homogeneous varieties.


\textbf{Mahler's method and Ku. Nishioka's theorem}

At the same time, K. Mahler introduces a method similar to the one of C. Siegel, which applies to (multivariate) functions, solutions of difference equations, rather than differential ones, of the form
$f\left(\Omega\bm{z}\right)=\frac{\sum_{l=0}^{m}a_{l}(\bm{z})f(\bm{z})^{l}}{\sum_{l=0}^{m}b_{l}(\bm{z})f(\bm{z})^{l}}$,
where $\bm{z}=(z_{1},\ldots,z_{n})$, the $a_{l}(\bm{z}),b_{l}(\bm{z})$ are polynomials with coefficients in a number field $\mathbb{K}$, that is a finite extension of $\mathbb{Q}$, and $\Omega=\left(r_{i,j}\right)\in\mathcal{M}_{n}(\mathbb{N})$ acts over $\mathbb{C}^{n}$ as
$\Omega\bm{z}=\left(\prod_{j=1}^{n}z_{j}^{r_{1,j}},\ldots,\prod_{j=1}^{n}z_{j}^{r_{n,j}}\right)$.

Moreover, K. Mahler requires that the biggest eigenvalue $\rho$ associated with $\Omega$ is real and that $1\leq m<\rho$.

The first main result obtained by the author guarantees, under some conditions, the transcendence of the values of a non-rational Mahler function at algebraic points.

%

However, contrary to $E$-functions, values of Mahler functions are not greatly prized numbers for applications. For example, we cannot obtain, a priori, $\pi$ as a value of such functions at an algebraic point. Yet, K. Mahler had two stimulating goals for his method. The first one concerned functions of the type of the Jacobi $\Theta$ function, as for example $\Theta(z_{1},z_{2})=\sum_{n=0}^{+\infty}z_{1}^{2n}z_{2}^{n^{2}}$. 

Indeed, $\Theta(\Omega(z_{1},z_{2}))=\frac{\Theta(z_{1},z_{2})-1}{z_{1}^{2}z_{2}}$, where $\Omega=\begin{pmatrix}
1 & 1\\
0 & 1
\end{pmatrix}$.

But the greatest eigenvalue of $\Omega$ is $\rho=1$, contrary to the assumptions of K. Mahler. 

The second target of K. Mahler was the Fourier expansion of the modular invariant $J$ attached to an elliptic curve, of which he conjectures the transcendence at algebraic points of the unit disc, except zero. Indeed, this function satisfies for every $q\in\mathbb{N}^{*}$ an equation of the form $P_{q}(J(z),J\left(z^{q}\right))=0$, where $P_{q}(X,Y)\in\mathbb{Z}[X,Y].$ K. Mahler mentioned but never studied these kind of polynomial equations, but Ku. Nishioka obtained transcendence results for values of solutions of such equations, under some assumptions. Unfortunately, these are not satisfied by $P_{q}$. Nevertheless, it is interesting to note that this conjecture from K. Mahler was proved in 1996 by K. Barré-Sirieix, G. Diaz, F. Gramain et G. Philibert using similar steps as the ones of Mahler's method. Thus, we can consider it as a fruitful source of inspiration. 

Mahler's method finally emerges in 1969, when W. Schwarz proves transcendental results for values of some functions, without mentioning nor exploiting the fact that they are Mahler functions. Then, K. Mahler publishes a paper in which he recalls his method and explain why the results of W. Schwarz are in fact corollaries of his own work. The author also suggests three open problems to generalise his approach. This gives rise to several works, from K. Kubota and J. Loxton and A. van der Poorten, who extend some results from K. Mahler and contribute to the popularisation of his method. Further work is made later by M. Amou, P-G. Becker, Ku. Nishioka, T. Töpfer, and many others. Through their impulse, the study moves to linear equations of the form
\begin{equation}
\label{eq_mahler_0}
P_{0}(z)f(z)+P_{1}(z)f(z^{d})+\cdots+P_{m}(z)f(z^{d^{m}})=0,
\end{equation}
where $P_{0}(z),\ldots, P_{m}(z)\in\mathbb{K}[z]$, $P_{m}(z)$ $\neq0$ and $\mathbb{K}$ is a number field.

Moreover, this new enthusiasm around Mahler's method is perpetuated by the links it shares with the dynamical research area of finite automata theory, especially spotlighted by M. Mendès-France. Indeed, every generating series of an automatic sequence is a Mahler function, in the sense of Equation \eqref{eq_mahler_0}.

The apotheosis of this first period of research is reached by Ku. Nishioka, who proves in 1990 the analogue of the Siegel-Shidlovskii's theorem for solutions of Mahler systems. Thereafter, these two results are refined into statements saying that, under the same respective assumptions, every homogeneous algebraic relation over $\Kbar$ between the numbers $\fialpha$ arises from the specialisation at $\alpha$ of a homogeneous algebraic relation over $\Kbarz$ between the functions $\fiz$ solutions of the differential or Mahler system concerned. For $E$-functions, this refinement is due to F. Beukers in 2006. For Mahler functions, it is obtained in 2017 by B. Adamczewski and C. Faverjon as a consequence of the work of P. Philippon. Recently, L. Nagy and T. Szamuely recovered this result from a generalization of works from Y. André.

\section*{Characteristic $p$}

 A big part of the classical theory of transcendence and algebraic independence deals with periods, in the sense of M. Kontsevich and D. Zagier. Examples of these remarkable numbers are $\pi$, values of logarithm at algebraic points, of the Riemann Zeta function at integral points, or powers of values of the Euler Gamma function at rational points. In characteristic $p$, there exists a framework in which we can build a parallel of this arithmetical context and define analogues of the complex exponential and logarithm functions, the Riemann Zeta function or the Euler Gamma function, along with associated values as $\pi$. This is the function fields setting, which shares profound analogies with the number fields one. Let us illustrate this fact by the following picture.
 
 \begin{figure}[H]
 	\begin{center}
 		\begin{tikzpicture}[scale=0.7]
 		
 		\node(Z) at (1,6) {$A=\mathbb{Z}$}; 
 		\node(Q) at (1,4.5) {$K=\mathbb{Q}$}; 
 		\node(Q_abs)[scale=0.75] at (3,4.5) {$|\frac{p}{q}|_{\infty}$ classical}; 
 		\node(Kd) at (1,3) {$\mathbb{K}$};
 		\node(Kdb) at (1,1.5) {$\overline{\mathbb{K}}$};
 		\node(Cd) at (1,0) {$C=\mathbb{C}$};
 		
 		\node(Rd) at (3,2.25) {$R=\mathbb{R}$};
 		
 		\node(car0) at (2,8) {Number fields};
 		\node(carp) at (8,8) {Function fields};
 		\node(carp') at (8,7.5) {of characteristic $p>0$};

 		\node(trait0) at (5,8) { }; 
 		\node(trait1) at (5,-0.5) { }; 
 		
 		\node(A) at (6.5,6) {$A=\mathbb{F}_{q}[T]$,};
 		\node(B) at (9,6) {$q=p^{r}$};
 		\node(K) at (6.5,4.5) {$K=\mathbb{F}_{q}(T)$}; 
 		\node(K_abs)[scale=0.75] at (8.65,4.85)	{$|\frac{P(T)}{Q(T)}|_{\infty}=$};
 		\node(K_abs)[scale=0.75] at (9.3,4.10)	{$\left(\frac{1}{q}\right)^{\deg_{T}(Q)-\deg_{T}(P)}$};
 		\node(Kd') at (6.5,3) {$\mathbb{K}$};
 		\node(Kdb') at (6.5,1.5) {$\overline{\mathbb{K}}$};
 		\node(C) at (6.5,0) {$C$};
 		
 		\node(R) at (8.9,2.8) {$R=\mathbb{F}_{q}\left(\left(\frac{1}{T}\right)\right)$};
 		\node(Rb) at (8.9,0.8) {$\overline{\mathbb{F}_{q}\left(\left(\frac{1}{T}\right)\right)}$};

 		\draw[-] (Z) -- (Q);
 		\draw[-] (Q) -- (Kd);
 		\draw[-] (Kd) -- (Kdb);
 		\draw[-] (Kdb) -- (Cd);
 		\draw[-] (Q) -- (Rd);
 		\draw[-] (Rd) -- (Cd);
 		
 		\draw[-] (trait0) -- (trait1);
 		
 		\draw[-] (A) -- (K);
 		\draw[-] (K) -- (Kd');
 		\draw[-] (Kd') -- (Kdb');
 		\draw[-] (Kdb') -- (C);
 		\draw[-] (K) -- (R);
 		\draw[-] (R) -- (Rb);
 		\draw[-] (Rb) -- (C);
 		
 		\node(f)[scale=0.75] at (0.5,3.7) {$<\infty$};
 		\node(fp)[scale=0.75] at (6,3.7) {$<\infty$};
 		
 		\node(c1)[scale=0.75] at (1.9,3.2) {c.};
 		\node(ca1)[scale=0.75] at (2.4,1) {a.c.};
 		\node(c2)[scale=0.75] at (7.4,3.6) {c.};
 		\node(c3)[scale=0.75] at (9.3,1.8) {a.c};
 		\node(ca2)[scale=0.75] at (7.2,0.45) {c.};

 		\end{tikzpicture}
 		\caption*{c. Completion for $|.|_{\infty}$. \\a.c. Algebraic closure.}
 	\end{center}
 \end{figure}

On the left, finite extensions $\mathbb{K}$ of $K$ are number fields. On the right, they are called function fields (of one variable). Periods are there replaced by remarkable numbers called periods of Drinfeld modules. They are any element of the free $A$-lattice $\Lambda\subseteq C$ of finite rank formed by the zeros of an exponential function defined in this setting. 

As an illustration, the analogue of the complex exponential function is called the Carlitz exponential. The associated lattice is $\Lambda=\varpi A$, with $\varpi=T(-T)^{1/(q-1)}\prod_{i=1}^{+\infty}\left(1-T^{1-q^{i}}\right)^{-1}$.\\ Consequently, the period $\varpi$ is an analogue of $2i\pi$ and $\pi_{q}:=\prod_{i=1}^{+\infty}\left(1-T^{1-q^{i}}\right)^{-1}$ is an analogue of $\pi$. 

A second illustration is the case of Drinfeld modules whose associated lattice has rank $2$. Indeed, they are analogues of elliptic curves, and we can attach to them periods, quasi-periods and a $j$-invariant. In 1998, M. Ably, L. Denis and F. Recher, obtain the transcendence of values of the Taylor expansion $J$ associated with $j$ for algebraic numbers $\alpha $ such that $0<|\alpha|<(1/q)^{1/(q-1)}$. This is the analogue in positive characteristic of the result due to K. Barré-Sirieix, G. Diaz, F. Gramain and G. Philibert in characteristic zero, mentioned above.

Beyond these similarities, there exist difficulties that specifically appear in the setting of function fields of positive characteristic. A significant example is the fact that the statement of K. Roth and consequently the one of the sub-space theorem of W. Schmidt is no longer satisfied. In spite of that, algebraic independence results in the function fields setting are often more developed than their analogues in the number fields framework, which are still conjectures. Let us now describe the methods developed in characteristic $p$ to obtain results of transcendence and algebraic independence of remarkable numbers in this setting. We will compare these statements to their analogues in characteristic zero.

\textbf{First results of L. Wade and M. Geijsel contribution}

The first main results dealing with transcendence of periods of Drinfeld modules are due to L. Wade in the 1940's. To begin with, the author proves that the values at non-zero algebraic points of the Carlitz exponential and logarithm are transcendental over $\overline{\FqTc}$. We recognize in the first statement the analogue of the classical C. Hermite's theorem. L. Wade also shows that the period $\xi$ he associates to the Carlitz exponential is transcendental. Then, J. Geijsel states that the set of the periods of the Carlitz exponential is $\mathbb{F}_{q}[T]\xi$. This extends the result of transcendence of L. Wade to all periods of the Carlitz exponential $e_{C}$. An other key result of this time is the analogue of the Gelfond-Schneider's theorem, obtained by L. Wade, and generalised by J. Geijsel. It states that, given a non-zero element $\alpha\in C$ such that $e_{C}(\alpha)\in\overline{\FqTc}$ and $\beta\in\overline{\FqTc}\setminus\FqTc$, the number $e_{C}(\alpha\beta)$ is transcendental over $\overline{\FqTc}$. 

\textbf{The automatic method}
\label{sec_meth_auto}

Around 1990, J-P. Allouche develops a new method to prove the transcendence of some numbers of the type $\alpha=\sum_{n=0}^{+\infty}a_{n}(1/T)^{n}\in \mathbb{F}_{q}((1/T))$. In particular, this strategy allows J-P. Allouche to recover the transcendence of $\pi_{q}$ first obtained by L. Wade. This approach is then generalised by V. Berthé, who recovers the transcendence over $\overline{\mathbb{F}_{q}(T)}$ of quotients $\zeta_{C}(s)/\pi_{q}^{s}$ for $s\in\{1,\ldots,q-2\}$, and the one of the elements $\log_{C}(P)/\pi_{q}^{s}$, for $s\in\{1,\ldots,q-2\}$, $q>3$ and certain rational fractions $P\in\mathbb{F}_{q}(T)$.

\begin{bandbox}[dcolor]
	The method of J-P. Allouche is based on two results. First, the one of G. Christol, which states that $\alpha$ is algebraic over $\overline{\mathbb{F}_{q}(T)}$ if and only if $(a_{n})_{n}$ is $q$-automatic. The other result is the one of S. Eilenberg which gives the following useful characterisation of automaticity. The sequence $(a_{n})_{n}$ is $k$-automatic if and only if its $k$-kernel $K_{k}(\underline{a})$ is finite, where $K_{k}(\underline{a})=\{\left(a_{k^{r}n+j}\right)_{n},r\geq 0, 0\leq j<k^{r}\}$.
\end{bandbox}

\textbf{The method of G. Anderson, W. D. Brownawell and M. Papanikolas}
\label{sec_methode_ABP}

A very powerful and fruitful method in the framework of function fields of positive characteristic is the one developed by G. Anderson, W. D. Brownawell and M. Papanikolas. It deals with $\tau$-difference systems, that is
\begin{equation}
\label{syst_ABP}
\begin{pmatrix}
\tau f_{1}(z) \\ \vdots \\ \tau f_{n}(z) 
\end{pmatrix} =A(z)\begin{pmatrix}
f_{1}(z) \\ \vdots \\ f_{n}(z)
\end{pmatrix},
\end{equation}
where each $f_{i}(z)$ admits a series expansion around zero, with coefficients in a function field of characteristic $p>0$ and where for $f(z)=\sum_{n=0}^{+\infty}a_{n}z^{n}$, with $a_{n}\in C$, we set $\tau f(z)=\sum_{n=0}^{+\infty}a_{n}^{1/q}z^{n}$.

This method is based on the $t$-motives theory introduced by G. Anderson around 1986, which generalises the notion of Drinfeld modules. This gives rise to analogues of classical results. For example, J. Yu uses this method to establish in 1997 the analogue of the well-known A. Baker's theorem. Besides, the three previous authors prove in 2004 the analogue of the refinements of the Siegel-Shidlovskii and Ku. Nishioka's theorems, for systems \eqref{syst_ABP}. 

In addition, the approach of the authors provides results which widely overcome the classical results of algebraic independence. Let us describe some of them. To begin with, L. Denis manipulates this theory to state in 1995 that $e=e_{C}(1)$ and $\varpi$ are algebraically independent over $\overline{\FqTc}$.

Moreover, C-Y. Chang and J. Yu consider $t$-motives to give in 2007 a complete description of all the algebraic relations over $\overline{\FqTc}$ between values at integers of the analogue in this setting of the Riemann function, called the Carlitz Zeta function.

Furthermore, M. Papanikolas uses this method to complete in 2008 the analogue of the A. Baker's theorem. More precisely, he states that, if $\alpha_{1},\ldots,\alpha_{n}\in\overline{\FqTc}$ are in the disc of convergence of the Carlitz logarithm $\log_{C}$ and if $\log_{C}(\alpha_{1}),\ldots,\log_{C}(\alpha_{n})$ are linearly independent over $\FqTc$, then these logarithms are algebraically independent over $\overline{\FqTc}$. 


A lot has also been found concerning algebraic independence of periods $\omega_{1},\omega_{2}$ and quasi-periods $\eta_{1},\eta_{2}$ attached to Drinfeld modules $\phi$ with lattice $\Lambda$ of rank $2$. The main result is due to C-Y. Chang and M. Papanikolas in 2012. The authors prove that $\omega_{1},\omega_{2},\eta_{1},\eta_{2}$ are algebraically independent, in the case where the characteristic $p$ is odd and $\phi$ is without complex multiplication. This means that we have $\{c\in C,c\Lambda\subseteq\Lambda\}=\FqTa$. The complex multiplication case, that is, $\FqTa\subsetneq\{c\in C,c\Lambda\subseteq\Lambda\}$, is different. Indeed, A. Thiery proves in 1992 that then, the transcendence degree over $\overline{\FqTc}$ of the family $\{\omega_{1},\omega_{2}, \eta_{1},\eta_{2}\}$ equals $2$.

\begin{bandbox}[dcolor]
	Concerning the result of C-Y. Chang and M. Papanikolas, in addition to the use of $t$-motives, their proof relies on the fundamental theorem of G. Anderson, W. D. Brownawell and M. Papanikolas mentioned above (which traces every homogeneous algebraic relation among values of solutions of $\tau$-difference systems back to a functional relation), and a strong result of M. Papanikolas, which turns the study of algebraic relations among periods and quasi-periods of Drinfeld modules $\phi$ into the study of a difference Galois group attached to a $\tau$-difference system associated to $\phi$.
\end{bandbox}

By comparison, in the classical case, we are very far from such results. Indeed, the algebraic independence over $\Qbar$ of $e$ and $\pi$ or the one of values of logarithms are conjectured but not proved. Concerning values of the Riemann $\zeta$ function, the only general result is the transcendence of its values at even integers. R. Apéry establishes in 1978 that $\zeta(3)$ is irrational. Since then, it has been proved by K. Ball and T. Rivoal in 2001 that an infinite number of values of $\zeta$ at odd integers is irrational. However, we still do not know how to find other irrational value apart from $\zeta(3)$. The closest result from this goal is due to W. Zudilin in 2004 and guarantees that at least one of the four numbers $\zeta(5), \zeta(7), \zeta(9), \zeta(11)$ is irrational.

In the same vein, our knowledge concerning elliptic curves is limited. Let us first mention T. Schneider who proves in 1957 the transcendence of every period associated to an elliptic curve. An other significant result is due to G. Chudnovsky in 1984 and gives the algebraic independence of $\frac{\pi}{\omega}$ and $\frac{\eta}{\omega}$, for every pair $(\omega,\eta)$ of period and quasi-period of an elliptic curve, when some invariants attached to its lattice are algebraic. Under this latter assumption, this result implies the existence of at least two algebraically independent numbers among generator periods and associated quasi-periods. But it seems that no more general result appears in this setting.

\textbf{Mahler's method in positive characteristic}

One of the advantages of Mahler's method rests in the fact that its validity does not depend on the characteristic involved. Moreover, a fruitful observation by L. Denis strongly motivates the use of this method in positive characteristic. It is the fact that, in the setting of function fields of positive characteristic, some periods, as $\pi_{q}$, arise from specialisations of Mahler functions at algebraic points. This makes Mahler functions quite natural in this framework, more than in the classical case in this sense. Through the innovative and prolific directions of the work of L. Denis, Mahler's method contributes since the 1990's to the first statements of algebraic independence that overcomes the classical case


But beyond the analogies between number fields and function fields of positive characteristic, some established results in the second setting do not admit any translation in the first one, because they concern objects that have no classical equivalents. For example, let us notice that the analogue $\pi_{q}\in \mathbb{F}_{q}((1/T))$ of $\pi$ depends on $q$. Thus, by changing $q$, we obtain several analogues of $\pi$. Then, L. Denis shows that, if $i_{1},\ldots,i_{n}$ is a strictly increasing sequence of non-zero natural numbers, then $\pi_{q^{i_{1}}},\ldots,\pi_{q^{i_{n}}}$ are algebraically independent over $\overline{\FqTc}$.  

Finally, using Mahler's method, G. Fernandes recently establishes the analogues of the Ku. Nishioka's theorem and its refinement in the function fields setting.

\section*{Related topics and perspectives}

\textbf{Algebraic independence measures}


\begin{bandbox}[dcolor]
Let us consider either one of the two settings introduced on page 3. An algebraic independence measure of numbers $\xi_{1},\ldots,\xi_{n}\in C$ is a lower bound of the form 
$|P(\xi_{1},\ldots,\xi_{n})|\geq\phi(\deg(P),H(P))$, where $\phi$ is a real function, which is satisfied for every polynomial $P\in\overline{\mathbb{K}}[\Xis]$ such that $P(\xi_{1},\ldots,\xi_{n})\neq 0$. Recall that $\deg(P)$ is the total degree of $P$ and $H(P)$ is the maximal of the absolute values of the coefficients of $P$.

\end{bandbox}


Let us restrict ourselves to values of $E$-functions, Mahler functions, and solutions of $\tau$-difference systems in positive characteristic. 

As far as we know, the third case did not generate general results on this topic. 


In the first two cases, in characteristic zero, such algebraic independence measures are due to A. Shidlovskii (using Siegel's method), and Ku. Nishioka (using Mahler's method) respectively, for algebraically independent functions. E. Zorin is now working on a general method, which, under some conditions, would in particular generalize the measures of Ku. Nishioka to algebraically dependent Mahler functions, in characteristic zero. Finally, in positive characteristic, only partial results are available. For example, P-G. Becker obtain in 1994 algebraic independence measure for Mahler functions with coefficients in $\mathbb{F}_{q}$. Furthermore, L. Denis uses Mahler's method in an interesting way in 2011 to get algebraic independence measures for some inhomogeneous Mahler equations of order $1$. A promising perspective of study could be to adapt the works of E. Zorin and L. Denis to obtain general algebraic independence measures for Mahler functions in positive characteristic.

\textbf{Galois theory}

 If we focus on methods developed by C. Siegel, Ku. Nishioka and  G. Anderson, W. D. Brownawell and M. Papanikolas, we note that the results and refinements obtained turn the problem of the algebraic independence of numbers into the one of the algebraic independence of functions. This question is at the heart of Galois theory of functional systems. Indeed, its main idea is to attach to the concerned system a linear algebraic group that reflects the algebraic relations between the solutions. This group is called the Galois group associated to the system. 
 
 Some authors apply the general ideas of this theory in the particular cases of $E$-functions, Mahler functions, or solutions to $\tau$-difference systems in positive characteristic. Concerning the first setting, the Hrushovski algorithm developed in 2002 computes the algebraic relations between solutions of differential systems. Very recently, S. Fischler and T. Rivoal adapt it (starting from the modified version by R. Feng) and deduce from Beuker's theorem an algorithm which computes a set of generators of the algebraic relations over $\Qbar$ between values of $E$-functions at algebraic points. For Mahler functions in characteristic zero, the work of J. Roques proves fundamental properties of the associated Galois group. The author also studies in detail the case of Mahler equations of order $2$. Finally, for $\tau$-difference systems in characteristic $p$, as mentioned above, a Galois theory is developed by M. Papanikolas. Besides, let us mention the Tannakian Galois theory, treated in particular by C. Hardouin. This approach is very general and especially includes the differential and $\tau$-difference cases. However, it seems that no Galois theory has been developed in the precise case of Mahler functions in positive characteristic. This could be a further perspective of study.
 
 To conclude, we shall stress that, even if the Galois theory is a powerful tool to deal with algebraic independence of functions that satisfy some functional equations, it is in general difficult to compute the associated Galois group. That is why, so far, explicit results of functional algebraic independence do not overcome equations of order $2$. The way is then widely open to fabulous discoveries !

\begin{wrapfigure}{l}{3.2cm}
\includegraphics[height=3.6cm, width=3.2cm]{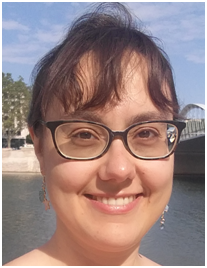} 
\end{wrapfigure} 
\textbf{{\large Gwladys Fernandes}}  \\~\\
Gwladys Fernandes is a Hadamard lecturer in mathematics at Versailles University. She recently defends her thesis on Mahler's method in positive characteristic. Her main research interests are in the theory of transcendence and algebraic independence and Galois theory for difference systems in positive characteristic. She is also curious about automata theory.

\end{twoblock}

%
%

\end{document}